\documentclass[11pt]{article}

\usepackage{amsmath,amsthm,amstext,amscd,amssymb,euscript,url}
\usepackage{mathrsfs}
\usepackage{calrsfs}
\usepackage{epsfig}
\usepackage[inline]{enumitem}
\usepackage[blocks]{authblk}
\usepackage{amsmath,amsthm,amstext,amscd,amssymb,euscript,url,mathrsfs,euscript,dsfont,calrsfs,mathtools,comment,relsize,xcolor}
\usepackage{hyperref}

\newcommand{\cov}{\text{cov}}
\newcommand{\Z}{\mathbb Z}
\newcommand{\R}{\mathbb R}
\newcommand{\N}{\mathbb N}

\newcommand{\E}{\mathbb E}

\renewcommand{\phi}{\varphi}

\newcommand{\La}{\ensuremath{\Lambda}}
\newcommand{\la}{\ensuremath{\Lambda}}

\newcommand{\loc}{\mathcal{L}}

\def\1{{\mathchoice {\rm 1\mskip-4mu l} {\rm 1\mskip-4mu l}
{\rm 1\mskip-4.5mu l} {\rm 1\mskip-5mu l}}}

\newtheorem{theorem}{{\small T}{\scriptsize HEOREM}}[section]
\newtheorem{corollary}{{\bf{\small C}{\scriptsize OROLLARY}}}[section]
\newtheorem{proposition}{{\bf{\small P}{\scriptsize ROPOSITION}}}[section]
\newtheorem{lemma}{{\bf{\small L}{\scriptsize EMMA}}}[section]
\newtheorem{remark}{{\bf{\small R}{\scriptsize EMARK}}}[section]
\newtheorem{definition}{{\bf{\small D}{\scriptsize EFINITION}}}[section]

\renewenvironment{proof}[1]
{\noindent{{\bf{\small{ P}{\scriptsize ROOF}}}.}\hspace{0.1cm} #1} {$\;\qed$\newline}

\newcommand{\beq}{\begin{eqnarray}}
\newcommand{\eeq}{\end{eqnarray}}

\newcommand{\ba}{\begin{align*}}
\newcommand{\ea}{\end{align*}}

\newcommand{\be}{\begin{equation}}
\newcommand{\ee}{\end{equation}}

\newcommand{\bl}{\begin{lemma}}
\newcommand{\el}{\end{lemma}}

\newcommand{\br}{\begin{remark}}
\newcommand{\er}{\end{remark}}

\newcommand{\bt}{\begin{theorem}}
\newcommand{\et}{\end{theorem}}

\newcommand{\bd}{\begin{definition}}
\newcommand{\ed}{\end{definition}}

\newcommand{\bp}{\begin{proposition}}
\newcommand{\ep}{\end{proposition}}

\newcommand{\bc}{\begin{corollary}}
\newcommand{\ec}{\end{corollary}}

\newcommand{\bpr}{\begin{proof}}
\newcommand{\epr}{\end{proof}}

\newcommand{\bi}{\begin{itemize}}
\newcommand{\ei}{\end{itemize}}

\newcommand{\ben}{\begin{enumerate}}
\newcommand{\een}{\end{enumerate}}

\newcommand{\caG}{{\mathcal G}}

\newcommand{\caR}{{\mathcal R}}

\newcommand{\caU}{{\mathcal U}}
\newcommand{\caV}{{\mathcal V}}

\newcommand{\IRW}{\text{\normalfont IRW}}

\newcommand{\SEP}{\text{\normalfont SEP}}

\newcommand{\RW}{\text{\normalfont RW}}

\begin{document}
\title{{\bf Particle systems with sources and sinks}}

\author[1]{Frank Redig
}
\author[2]{Ellen Saada
}

\affil[1]{Institute of Applied Mathematics, Delft University of Technology, Delft,
The Netherlands}
\affil[2]{CNRS, UMR 8145, Laboratoire MAP5, Universit\'{e} Paris Cit\'{e},
Paris, France}

\date{Dated: \today}

\maketitle

\begin{abstract}
Local perturbations in conservative particle systems can have a
non-local influence on the stationary measure. To capture this phenomenon,
we analyze 
two toy models.
We study the symmetric  exclusion process  on a countable set of sites
$V$ with a source at a given point (called the origin),
starting from a Bernoulli product measure with density $\rho$.
We prove that when the underlying random walk  on $V$  is recurrent, then
the system evolves towards full occupation, whereas in the transient case
we obtain a limiting distribution which is not product and
has long-range correlations.
For independent random walkers  on $V$,  we analyze the same problem,
starting from a Poissonian
measure. Via intertwining with a system
of ODE's, we prove that the distribution is  Poissonian
at all later times $t>0$, and that
the system ``explodes'' in the limit $t\to\infty$ if and only if the underlying
random walk is recurrent. In the transient case, the limiting density is
a simple function of the Green's function  of the random walk.
\end{abstract}

\section{Introduction}
Introducing local perturbations in conservative particle systems can have
a drastic, i.e.,
non-local influence on the stationary measure, which can change from
 a product form
to a measure with long-range correlations.
This has already been studied in \cite{maes, maesredig}
where anisotropic perturbations
of the symmetric  exclusion process
were shown to exhibit long-range correlations,
using a formal series
expansion method. See also the more recent work \cite{muka},
where an exclusion process with a driven bond is studied,
and \cite{krap} where the growth of the total number of particles
in a symmetric exclusion process with source is studied.
The abelian sandpile model \cite{dhar} is another well-known important example
showing that a conservative diffusive dynamics combined with sources leads
to a self-organized critical state,  i.e., a stationary measure which has
power law decay of correlations.

 To capture this phenomenon, we analyze in this paper two toy models.
We first consider the symmetric  exclusion process
with a source at the origin,
in the setting of an infinite graph, and start it from a Bernoulli
product measure of constant density,
which is the stationary probability measure for the dynamics  without source.
In case the underlying random walk  of the exclusion process is recurrent,
we show that the system becomes fully occupied, whereas in the transient case,
a limiting measure is obtained.
This measure is the microscopic analogue of the solution of the Poisson equation,
 and it is shown to have long-range correlations.
We then study a similar setting where the source at the origin is replaced
by the combination of a source and a sink, which can be thought of as
 a coupling to a reservoir. For the density corresponding to the reservoir density
 the stationary measure is a product measure, but for other initial densities,
 in the transient case another limiting measure (with identical limiting density)
 with long-range correlations is obtained.  These results show 
 that a reservoir coupled to an infinite system directly has very different 
 effects when compared to the standard setting of boundary reservoirs which 
 produce a non-equilibrium steady state. The main technique of proof is a 
 combination of duality with the Feynman-Kac formula.
 These techniques can be also used in a more general setting of a finite 
 number of sources and sinks. 

Finally, we consider the same setting in the case of independent random walkers,
where we show that starting from a homogeneous product of
Poissonian distributions, the distribution at any later time is still
a  Poisson product measure, with a density which diverges in the recurrent case,
and has a limit in the transient case. The main tool for the proof is
a  new intertwining relation 
between the independent walkers with a source and
a deterministic system of coupled linear differential equations
which can be solved explicitly.

The rest of our paper is organized as follows. In section \ref{sec:settings}
 we introduce the exclusion process with source (and possibly sink) at the origin,
 and explain the main questions. In section \ref{sec:inv-for-source}
 we state and prove the main result, that is, about invariant measures
 for the exclusion processes with source or with source and sink.
 In section \ref{sec:further} we prove more on these measures 
 for the process with a source,
 namely negative correlations and computation of covariances.
  In section \ref{sec:irw} we consider the case of independent random walkers.

\section{Setting, notations and definitions}\label{sec:settings}
Let $V$ denote a countable set of vertices  and
$(p(x,y), x,y\in V)$
 irreducible and symmetric random walk transition rates on $V$,
i.e., for every $x,y\in V$, $p(x,y)=p(y,x)\geq 0$,  and
 there exists $n\in \N$ such that $p^{(n)}(x,y)>0$,
where $p^{(n)}(x,y)$ denotes the corresponding $n$-step
transition rate, i.e.,
\[
p^{(n)}(x,y)=\sum_{z_1,\ldots,z_{n-1}\in V} p(x,z_1)\ldots p(z_{n-1},y),
\]
Note that we do not assume a priori that $(p(x,y), x,y\in V)$ is a probability
transition.
In order to avoid existence problems, we will assume that
\[\sup_{x\in V} \sum_{y\in V} p(x,y)<\infty.
\]
\subsection{Single particle dynamics}\label{subsec:single}
 We denote by $\{X_t:t\geq 0\}$ a continuous-time random walk moving with rate
 $(p(x,y): x,y\in V)$ over the (oriented) edge $xy$, i.e.,
 the random walk on $V$ with generator
\be\label{Agen}
 L_{\RW}f(x)  = \sum_{y\in V} p(x,y) (f(y)-f(x)),
\ee
for a function $f:V\to \R$.
We denote by $\E^{\RW}_x$ expectation for this random walk starting at $X_0=x$.

\subsection{Exclusion process}\label{subsec:exclusion}
Next, we consider the symmetric exclusion process  on $V$  based on
$(p(x,y): x,y\in V)$,  which intuitively speaking consists
of independent walkers
moving each one according to the generator  $L_{\RW}$
and subject to the restriction
that at any instant of time multiple occupancies are forbidden, i.e., all
jumps leading to more than one particle per site are forbidden.

We thus consider configurations of particles with at most one particle
per site and denote the corresponding configuration space by $\Omega= \{ 0,1\}^V$.
Elements of $\Omega$ are denoted by $\eta, \xi,\zeta$, and for
$\eta\in \Omega$ we denote
by $\eta_x$  the occupation at vertex $x\in V$, i.e, $\eta_x=1$
(resp.\ $\eta_x=0$) means $x$ is occupied (resp.\ empty).

The exclusion process based on $(p(x,y): x,y\in V)$ is then defined as
the unique Markov process
$\{\eta(t):t\geq 0\}$  on $\{0,1\}^V$ with generator given by
\be\label{exgen}
 L_{\SEP} f(\eta)
= \sum_{x,y\in V} p(x,y) \eta_x (1-\eta_y) (f(\eta^{x,y})- f(\eta)),
\ee
for a local function $f: \Omega\to\R$ (i.e., depending
on a finite number of coordinates $\eta_i, i\in V$),
where $\eta^{x,y}$ is the configuration obtained from $\eta$ by removing
a particle at $x$ and putting it at $y$,  that is, for $\eta_x=1,\eta_y=0$
\be\label{def:eta_xy}
(\eta^{x,y})_z =
\begin{cases}
\eta_x -1 &\hbox{ if }  z=x\cr
\eta_y +1 &\hbox{ if }  z=y\cr
\eta_z &\hbox{ otherwise}.\cr
\end{cases}
\ee
Otherwise, we put $\eta^{x,y}=\eta$.
The existence of the process with generator \eqref{exgen}
is proved in \cite[Chapter VIII]{ligg}, that also contains
 many properties of the process,
including self-duality  (see next subsection).
The interpretation of the generator is
that particles  move in continuous time according
to the hopping rates $(p(x,y): x,y\in V)$
but jumps to already occupied sites are suppressed.

If in the exclusion process we start from a finite number of particles,
initially located at sites in a set
$A=\{x_1, \ldots, x_n\}\subset V$ then at any later time,
the set of occupied vertices
 is a finite set $A(t)=\{ x_1(t), \ldots, x_n(t)\}$.
 For $n=2$ we write $A=\{x,y\},\,A(t)=\{ x(t), y(t)\}$.
 Thus we equivalently  consider the exclusion process
as taking values in the set of finite subsets of $V$ and (with an abuse of notation)
we may write its generator also as
\be\label{exgen-bis}
 L_{\SEP} f(A) =\sum_{x,y\in V}p(x,y)I (x\in A,y\notin A)( f(A^{x,y}) - f(A)),
\ee
 for a local function $f$,
where $I$ denotes  the indicator function and
$A^{x,y}$ is the set obtained by removing $x$ from $A$ and adding $y$ to $A$.
We denote by $\E_A^{\SEP}$ the  expectation in the process $\{A(t),t\geq 0\}$,
by $\E^{\SEP}_{x_1,\ldots, x_n}$ expectation, by $\mathcal{V}_n$ the generator
and by $V_n(t)$ the semi-group
in the corresponding labeled process $(x_1(t), \ldots, x_n(t))$,
where we choose the initial labels, which are preserved in the course of time.
We thus write, for a local function $f$,
\[
V_n(t) f(x_1, \ldots, x_n)=\E_{x_1,\ldots, x_n}^{\SEP} f(x_1(t), \ldots, x_n(t)).
\]
A function of the configuration $\eta(t)$ can then alternatively be viewed
as a symmetric function of $x_1(t), \ldots, x_n(t)$.

For  $n$ independent random walkers
(each one with generator \eqref{Agen})
 initially located on $x_1, \ldots, x_n$, 
we denote respectively by $\mathcal{U}_n$, $U_n(t)$,
$\E^{\IRW}_{x_1,\ldots, x_n}$ the corresponding generator, semi-group and
expectation,  and by $X_{1,t}, \ldots, X_{n,t}$
 their locations at time $t>0$
 (for $n=2$ we write  $X_t,Y_t$).  More generally,
the independent random walkers process is then a process on $\N^V$
with formal generator
\be\label{ind}
L_{\IRW} f(\eta) =\sum_{x,y\in V} p(x,y)\left( \eta_x (f(\eta^{x,y})- f(\eta))
 + \eta_y (f(\eta^{y,x})- f(\eta))\right),
\ee
working on local functions $f: \N^V\to\R$. We denote by $\E^{\IRW}_\eta$ the
corresponding expectation, when
starting from the configuration $\eta$.

\subsection{Self-duality for exclusion and random walkers}\label{subsec:duality}
An important property of the (symmetric) exclusion process is self-duality
(see \cite[Chapter VIII]{ligg}), which is formulated as follows.
Let $\xi\in\Omega$ denote a finite configuration, i.e., $\sum_x \xi_x<\infty$,
and,  for $\eta\in\Omega$, define
\be\label{def:dual-fct-1}
D(\xi, \eta)= I(\xi\leq \eta),
\ee
where
$\xi\leq \eta$ refers to coordinate-wise order,  i.e.,
if there is a particle in $\xi$ at a site $x$,
 then there must also be a particle
in $\eta$ at $x$.
Then we have, for any finite configuration $\xi$,
for any $\eta\in\Omega$ and $t>0$
\be\label{sdualex}
 \E_\eta^{\SEP} D(\xi, \eta(t))= \E_\xi^{\SEP} D(\xi(t), \eta).
\ee
Let us denote by $\xi=e_x$ the configuration with a single particle
at $x$ and no particles elsewhere, then,
under the exclusion process $\{\xi(t), t\geq 0\}=  \{ e_{X_t}:t\geq 0\}$,
where $X_t$ is the random walk
with generator   $L_{\RW}$ of \eqref{Agen},
 the duality relation \eqref{sdualex}  reads
\be\label{sdualex2}
\E_\eta^{\SEP} \eta_x(t)= \E_x^{\RW} (\eta_{X_t}).
\ee
The symmetric exclusion process has as reversible
 probability measures
homogeneous Bernoulli product measures.
We denote them by
\be\label{def:nurho}
\nu_\rho, \quad\rho\in[0,1], \quad \text { with }
\nu_\rho(\eta_x=1)=\rho, \,\text { for all } x\in V.
\ee
Similarly, there is a self-duality relation
for independent random walkers, which reads as follows.
Define the polynomials
\[
D_{\IRW}(\xi, \eta)= I(\xi\leq \eta) \prod_{i\in V} \frac{\eta_i!}{(\eta_i-\xi_i)!},
\]
where $\xi\in \N^V$ is a finite configuration of
the independent random walkers process,
and $\eta\in\N^V$.
Then we have
\be\label{indepduality}
\E^{\IRW}_\eta (D_{\IRW}(\xi, \eta(t)))= \E^{\IRW}_\xi (D_{\IRW}(\xi(t), \eta)).
\ee
 For the proof of this well-known self-duality relation,
 we refer e.g. to \cite{dmp}.
\subsection{The model with a source
(or with source and sink)}\label{subsec:with-source}
To define the process with a source (or a source and a sink),
we fix a vertex $0\in V$, that we call the origin,
and call  \textit{symmetric exclusion process with source at $0$},
\textit{of intensity} $\lambda\geq 0$ the process with generator
\be\label{sourcexgen}
 L_{\SEP,0} f(\eta)
= \sum_{x,y} p(x,y) \eta_x (1-\eta_y) (f(\eta^{x,y})- f(\eta))
+ \lambda (1-\eta_0) (f(\eta^0)-f(\eta)),
\ee
for a local function $f$,
where $\eta^0$ is the configuration obtained from $\eta$ by
flipping the occupation variable at the origin,
 that is,
\be\label{def:etahat0}
(\eta^{0})_z =
\begin{cases}
1-\eta_0 &\hbox{ if }  z=0\cr
\eta_z &\hbox{ otherwise}.\cr
\end{cases}
\ee
In other words: particles move according to the symmetric exclusion process,
and whenever the origin is empty at rate $\lambda$ a particle is added.
 We denote by $S_0(t)$ the semi-group of this process.
We start the process with source at $0$ from $\nu_\rho$
 (defined in \eqref{def:nurho}),
 we denote by $\nu_{\rho,0} (t):=\nu_\rho S_0(t)$
the measure at time $t>0$.

Note that
we use the sub-index 0 to refer that the process has a source at the origin;
 later on,
 we will use the sub-index 1 for the process with a source and a sink
 at the origin. 

We then study the following two questions.
\ben
\item When is  $\lim_{t\to\infty} \nu_{\rho,0} (t)$
equal to $\delta_{\underline 1}$,
 the  Dirac measure concentrated on the
fully occupied configuration
 (that is, such that $\underline 1(x)=1$ for any site $x\in V$)?
\item If $\lim_{t\to\infty} \nu_{\rho,0} (t)\not= \delta_{\underline 1}$,
what is the limiting measure? What is its density,
and are there non-trivial correlations?
\een
The same questions will also be asked for a model with a source
and a sink at $0$, i.e.,
the process with generator
\beq\label{sourcsinkexgen}
 L_{\SEP,1} f(\eta)
 &=& \sum_{x,y} p(x,y) \eta_x (1-\eta_y) (f(\eta^{x,y})- f(\eta)) \nonumber\\
&+ & \lambda (1-\eta_0) (f(\eta^0)-f(\eta)) +  \mu\eta_0 (f(\eta^0)-f(\eta)),
\eeq
 for an intensity $\mu>0$, for a local function $f$.
 We denote by $S_1(t)$ the semi-group of this process.
We start the process with source and sink at $0$ from $\nu_\rho$,
 we denote by $\nu_{\rho,1} (t):=\nu_\rho S_1(t)$
the measure at time $t>0$.
Notice that for this process the Bernoulli product measure with density
\be\label{def:rhoR}\rho_R:= \frac{\lambda}{\lambda+\mu}\ee
 is reversible, but this fact
does not imply that this measure
is the only invariant probability measure, i.e., the two questions asked
for the models with source can be asked for the model with source and sink as well.
The idea is that the source and sink site corresponds to a ``reservoir''
and that in the transient case the system can ``miss the reservoir'', and
therefore converge to a limiting density different from the  density
$\rho_R$ imposed by the reservoir.

\section{Invariant measures}\label{sec:inv-for-source}
We consider the process with a source, that is,
with generator \eqref{sourcexgen}, 
 then the process with a source and a sink, 
that is with generator \eqref{sourcsinkexgen}, and in both cases 
 we denote by
$\E_\eta$ expectation in this process starting from configuration $\eta\in \Omega$.
We denote by $\E_{\nu_\rho}= \int \E_\eta \nu_\rho(d\eta)$ expectation
starting from an initial configuration $\eta$ which is $\nu_\rho$ distributed.
First we show the existence of $\lim_{t\to\infty} \nu_\rho S_0(t)=:\mu_{\rho,0}$.
This is proved with the help of a dual process, where particles
are moving according to an exclusion process with a sink at the origin.
 Then we look for 
$\lim_{t\to\infty} \nu_\rho S_1(t)=:\mu_{\rho,1}$. 
\subsection{Two duals  of the model with source
and convergence to an invariant measure 
for both models, with a source or with a source 
and a sink}\label{subsec:dual-and-conv}
 In this subsection we introduce two alternative duality relations
 for the process with a source,
analogous to the ones introduced in \cite[Chapter III]{ligg} for spin systems.
For the first one, that is a ``killed random walkers dual'',
denote, for $A\subset V$ a finite set,
\be\label{bimbam}
H(A,\eta)= \prod_{x\in A} \eta_x,
\ee
with the convention $ H(\emptyset, \eta)= 1$.
Then, we compute
\be\label{bobbi}
\lambda(1-\eta_0)\left(H(A,\eta^0)- H(A,\eta)\right)
= \lambda I (0\in A)( H(A\setminus \{0\},\eta)- H(A,\eta))
\ee
where we used
\be\label{bimbam-with-0}
H(A, \eta^0)=
\begin{cases}
H(A, \eta) &\hbox{ if } 0\not\in A \cr
\left(\prod_{x\not= 0, x\in A} \eta_x\right) (1-\eta_0) &\hbox{ if } 0\in A.\cr
\end{cases}
\ee
By combining \eqref{bimbam}--\eqref{bimbam-with-0} with the self-duality
relation
\eqref{sdualex}  of the symmetric exclusion process we find
\be\label{adual}
L_{\SEP,0} H(A, \eta) = \overline{\loc} H(A,\eta),
\ee
where the generator $\overline{\loc}$ works on the $A$-variable
(i.e., on finite subsets of $V$)  and is given by
\beq\nonumber
\overline{\loc} f(A)
&=& \sum_{x,y\in V}p(x,y)I (x\in A,y\notin A)( f(A^{x,y}) - f(A))\\\label{adualloc}
&&+ \lambda I( 0\in A) (f(A\setminus \{0\})- f(A)),
\eeq
 for a local function $f$.
The dual process  $\{ \overline{A}(t), t\geq  0\}$
with generator $\overline{\loc}$ is a process
taking values in the set of finite subsets of $V$, and can be described
as follows: particles initially located at the sites of $A$ perform
the symmetric exclusion process starting from $A$ (i.e., there are
particles at the sites of $A$ and no particles elsewhere),
and are killed with rate $\lambda$ when they are at the origin.
Let us denote $\E^{dual}_A$ expectation in this process starting from $A$,
and further denote $p_t^\lambda(A,B)$ the transition probability in this
dual process to go from $A$ to $B$ in time $t$; hence note that here,
as each time we introduce a transition $p_t(.,.)$, either for exclusion
 or for random walker, we assume that $(p(x,y), x,y\in V)$
 is a probability transition.\par\smallskip

As a consequence, we have the following.
\bt\label{dualitythm}
Let $A\subset V$  be a finite set.
Then
we have, for every $t>0$
\be\label{murhoha}
\int H(A,\eta) \nu_{\rho,0}(t) (d\eta)
= \sum_{B\subset V \atop |B|\leq |A|, B\not=\emptyset } p_t^\lambda (A, B) \rho^{|B|}
 + p_t^\lambda(A,\emptyset) =
 \E^{dual}_A(\rho^{|\overline{A}(t)|}),
\ee
and as $t\to\infty$ the limiting measure
  $\mu_{\rho,0} = \lim_{t\to\infty} \nu_{\rho,0}(t) $
exists, is invariant and satisfies
\be\label{intHA}
 \int H(A,\eta) \mu_{\rho,0}(d\eta) =  \E^{dual}_A(\rho^{|\overline{A}_\infty|}),
\ee
 where $|\overline{A}_\infty|=\lim_{t\to\infty}|\overline{A}(t)|$.
\et
\bpr
By the duality relation \eqref{adual} we obtain
\be\label{expHAt}
\E_\eta H(A,\eta(t))=  \E^{dual}_A H(\overline{A}(t), \eta).
\ee
Integrating this relation over $\nu_\rho$ (in the $\eta$-variable)
yields \eqref{murhoha}. The limit $t\to\infty$ is well-defined
because the cardinality of $\overline{A}(t)$  is non-increasing in $t$.
The invariance of $\mu_{\rho,0}$ follows because it is equal to the limit
$ \lim_{t\to\infty} \nu_{\rho,0}(t)$
(by \cite[Chapter I, Proposition 1.8]{ligg}).
\epr

There is also an alternative duality relation with a Feynman-Kac term, that
 can be obtained as follows.
If we define, for $A\subset V$ a finite set,
\be\label{def:tildeH}
\widetilde{H}(A, \eta)=\prod_{x\in A} (1-\eta_x),
\ee
 then we find
\[
L_{\SEP,0}\widetilde{H}(A, \cdot)(\eta)= \caR\widetilde{ H}(\cdot, \eta)(A),
\]
where the operator $\caR$ is in Schr\"{o}dinger operator form and given by
\[
\caR f(A)= \sum_{x,y\in V} I (x\in A,y\notin A)p(x,y)( f(A^{x,y}) - f(A))
 - \lambda I( 0\in A)  f(A),
\]
 for a local function $f$.
As a consequence of the Feynman-Kac formula we then have
\be\label{fedual}
\E_\eta \widetilde{H}(A, \eta(t))= \E_A^{\SEP}
\left( e^{-\lambda \int_0^t I(0\in A(s)) ds}
\widetilde{H}(A(t), \eta)\right).
\ee
Upon integrating \eqref{fedual} over the Bernoulli measure $\nu_\rho$ gives
\[
\int\widetilde{H}(A, \eta) \nu_{\rho,0}(t)(d\eta)
= (1-\rho)^{|A|} \E_A^{\SEP}
\left( e^{-\lambda \int_0^t  I(0\in A(s)) ds}\right).
\]

Then we have the following result on the limiting measure $\mu_{\rho,0}$
for the dynamics with a source, as well as on $\mu_{\rho,1}$
for the dynamics with source and sink.
\bt\label{sourcethm}
\bi
\item[I)] {\bf \em The model with source.}
Let $\{\eta(t), t\geq 0\}$ denote the process with generator
\eqref{sourcexgen}.
\ben
\item If $(p(x,y): x,y\in V)$ is recurrent then
$\mu_{\rho,0}=\delta_{\underline 1}$.
Moreover for any
configuration $\eta\in \Omega$, we have
$\lim_{t\to\infty} \delta_\eta S_0(t)=\delta_{\underline 1}$,
 where  $\delta_\eta$
denotes the Dirac measure on configuration $\eta$.
As a consequence in that case $\delta_{\underline 1}$
is the unique invariant  probability  measure.
\item If $(p(x,y): x,y\in V)$ is transient then
$\lim_{t\to\infty}\nu_{\rho,0} (t)= \mu_{\rho,0}$ with the following properties
\bi
\item[a)] Limiting density:
\be\label{dens}
\lim_{t\to\infty} \E_{\nu_\rho} (1-\eta_x(t))
= (1-\rho) \E^{\RW}_x e^{-\lambda \int_0^\infty \delta_{X_s,0}\ ds},
\ee
where $\delta_{\cdot,\cdot}$ denotes the Kronecker symbol.
\item[b)] Covariances:
\beq\label{covmu}
\cov_{\mu_{\rho,0}} (\eta_x, \eta_y)
&=& (1-\rho)^2 \left(\E_{x,y}^{\SEP}
e^{-\lambda \int_0^\infty (\delta_{x(s),0}+\delta_{y(s),0})\ ds}\right.\cr
 && \qquad\left.-\E_{x,y}^{\IRW}
 e^{-\lambda \int_0^\infty (\delta_{X_s,0}+\delta_{Y_s,0})\ ds}\right).
\eeq
\ei
\een
\item[II)] {\bf \em The model with source and sink.}
Let $\{\eta(t), t\geq 0\}$ denote the process with generator
\eqref{sourcsinkexgen}.
\ben
\item If $(p(x,y): x,y\in V)$ is recurrent then
$\lim_{t\to\infty}\nu_{\rho,1} (t)=  \nu_{\rho_R}$
(recall \eqref{def:rhoR}).  
The same holds for any initial configuration, i.e.,
$\lim_{t\to\infty}\delta_\eta (t)= \nu_{\rho_R}$ 
for every $\eta\in\Omega$.
\item If $(p(x,y): x,y\in V)$ is transient then
 $\lim_{t\to\infty}\nu_{\rho,1} (t)=\mu_{\rho,1}$  with the following properties
\bi
\item[a)] Limiting density:
\be\label{denstwo}
\lim_{t\to\infty} \E_{\nu_\rho}( \eta_x(t))
= \rho_R
+ \left(\rho-\rho_R \right)
\E^{\RW}_x e^{-(\lambda+\mu) \int_0^\infty \delta_{X_s,0}\ ds}.
\ee
\item[b)] We have the following formula for the covariances
\beq\nonumber
\cov_{\mu_{\rho,1}} (\eta_x, \eta_y)
&=& (\rho-\rho_R)^2 \left(\E_{x,y}^{\SEP} e^{-\lambda
\int_0^\infty (\delta_{x(s),0}+\delta_{y(s),0})\ ds}\right.\\\label{covmutwo}
&&\left.-\E_{x,y}^{\IRW}
e^{-\lambda \int_0^\infty (\delta_{X_s,0}+\delta_{Y_s,0})\ ds}\right).
\eeq
\ei
\een
\ei
\et
\bpr  \textbf{Part I, case 1.}
Let us start by the computation starting from the generator \eqref{sourcexgen}
for $x\in V$,
\be\label{bubu}
L_{\SEP,0} (1-\eta_x)= \sum_{y\in V} p(x,y) ((1-\eta_y)-(1-\eta_x))
-\lambda \delta_{x,0} (1-\eta_x).
\ee
Denoting
\be\label{not:psi}
\psi(x,\eta)=(1-\eta_x), 
\ee
 we can write \eqref{bubu} in the form 
\be\label{bibi}
L_{\SEP,0} \psi(x,\cdot) 
=  \sum_{y\in V} p(x,y)
(\psi(y,\cdot)-\psi(x,\cdot))  -\lambda \delta_{x,0} \psi(x,\cdot).
\ee 
Let us denote by  $\mathcal{A}_\lambda$ 
 the operator  (with parameter $\lambda$) 
 working on the $x$-variable as
\be\label{not:Ala}
 \mathcal{A}_\lambda  \phi(x)
=  \sum_{y\in V} p(x,y)(\phi(y)-\phi(x))  -\lambda \delta_{x,0} \phi(x),
\ee
which is the sum of the random walk generator $L_{\RW}$
and a multiplication operator with the ``potential'' $-\lambda \delta_{x,0}$.
As a consequence, by the Feynman-Kac formula we obtain
\be\label{femi}
e^{t \mathcal{A}_\lambda } \phi(x)
= \E_x^{\RW}\left(e^{-\lambda \int_0^t \delta_{X_s,0}\ ds}
\phi(X_t)\right).
\ee
This identity, combined with \eqref{bibi} gives
\be\label{cala}
\E_\eta(1-\eta_x(t))
= e^{t \mathcal{A}_\lambda } \psi(\cdot, \eta) (x)
= \E_x^{\RW}\left(e^{-\lambda \int_0^t \delta_{X_s,0}\ ds}(1-\eta_{X_t}) \right).
\ee
Integrating over $\nu_\rho$ gives
\be\label{cala-int}
\int \nu_\rho (d\eta) \E_\eta(1-\eta_x(t))
=(1-\rho) \E_x^{\RW}\left(e^{-\lambda \int_0^t \delta_{X_s,0}\ ds}\right).
\ee
Taking the limit $t\to\infty$ gives statements 1 (because the integral
on the right hand side of \eqref{cala-int} goes to $+\infty$),
 and  2a) of Part I of Theorem \ref{sourcethm}. \\

{\bf Part II, case 1}.
To prove the case with source and sink, start again from
the computation of $L_{\SEP,1}\eta_x$ for the generator
\eqref{sourcsinkexgen} and $x\in V$;
we get,
\be\label{bibibum}
 L_{\SEP,1} \eta_x  =  \sum_{y\in V} p(x,y) (\eta_y-\eta_x)
-(\lambda+\mu) \delta_{x,0} \eta_x + \lambda\delta_{x,0}.
\ee
Let us consider the equation  (recall the notation \eqref{not:Ala}) 
\be\label{koko}
\frac{d}{dt} \phi(x,t)=  
\mathcal{A}_{\lambda+\mu}  \phi(\cdot, t) (x) + \lambda \delta_{x,0},
\ee
Then by the variation of
constants method, we find the solution
\be\label{kokoriko}
\phi(x,t)
 = e^{t \mathcal{A}_{\lambda+\mu} }\phi(x,0)
 + \int_0^t e^{(t-s) \mathcal{A}_{\lambda+\mu} } \delta_{x,0}\ ds,
\ee
where $ \mathcal{A}_{\lambda+\mu} $
as well as $e^{t \mathcal{A}_{\lambda+\mu} }$
work on the $x$-variable.
Because the semigroup $e^{t \mathcal{A}_{\lambda+\mu} }$
can be computed using the Feynman-Kac formula we obtain, combining
\eqref{bibibum} with \eqref{kokoriko},
 using the analogue of \eqref{femi}, 
\beq
&&\E_\eta (\eta_x(t))
\nonumber\\ &=&  \E_x^{\RW}
\left(e^{-(\lambda+\mu) \int_0^t \delta_{X_s,0}\ ds}\eta_{X_t} \right)
+
\lambda\int_0^t \E_x^{\RW}
\left(e^{-(\lambda+\mu) \int_0^s \delta_{X_r,0}\ dr}\delta_{X_s,0} \right)\ ds
\nonumber\\
&=& \nonumber
\E_x^{\RW}\left(e^{-(\lambda+\mu) \int_0^t \delta_{X_s,0}\ ds}\eta_{X_t} \right)\ ds
\\\nonumber&&+
\frac{\lambda}{\lambda+\mu}\int_0^t
\E_x^{\RW}\left(e^{-(\lambda+\mu)
\int_0^s \delta_{X_r,0}\ dr} (\lambda+\mu)\delta_{X_s,0} \right) \ ds
\nonumber\\
&=&
\E_x^{\RW}\left(e^{-(\lambda+\mu) \int_0^t \delta_{X_s,0}\ ds}\eta_{X_t} \right)\ ds
+
 \rho_R  \int_0^t (-1)\frac{d}{ds}
 (e^{s \mathcal{A}_{\lambda+\mu} } 1) (x)\ ds
\nonumber\\
&=&\nonumber
\E_x^{\RW}\left(e^{-(\lambda+\mu) \int_0^t \delta_{X_s,0}\ ds}\eta_{X_t} \right)\ ds
\\&&+
 \rho_R  \left( 1-
\E_x^{\RW}\left(e^{-(\lambda+\mu) \int_0^t \delta_{X_s,0}\ ds}\right)\right),
\eeq
which yields items 1 and 2a) of Part II of Theorem \ref{sourcethm}.\\

 {\bf Part I, case 2}.
We now focus on the transient case and prove statement 2b)
of part I of Theorem \ref{sourcethm}. We compute the expectation
\[
\E_\eta (\psi(x,\eta(t))\psi(y, \eta(t))),
\]
for  $x,y\in V$, 
where we remind the reader  the notation \eqref{not:psi}.
 First, a generator computation yields, using self-duality
of the symmetric exclusion process (see \eqref{sdualex},
\eqref{sdualex2}), for $x,y\in V$,
\be\label{popi}
 L_{\SEP,0}
\psi(x,\cdot)\psi(y,\cdot)=
\mathcal{V}_2 \psi(x,\cdot)\psi(y,\cdot)
 -\lambda (\delta_{x,0}+ \delta_{y,0}) \psi(x,\cdot)\psi(y,\cdot),
\ee 
where $\mathcal{V}_2$ denotes the generator of
two exclusion particles
initially starting from $x,y$ (cf. subsection \ref{subsec:exclusion}).
Using once more the Feynman-Kac formula, this leads to
\be\label{feika}
\E_\eta (\psi(x,\eta(t))\psi(y, \eta(t)))
= \E^{\SEP}_{x,y}\left( e^{-\lambda\int_0^t
(\delta_{x(s),0}+ \delta_{y(s),0})\ ds }\psi(x(t), \eta)\psi(y(t),\eta)\right).
\ee
Integrating this equality over $\nu_\rho$ (in the $\eta$-variable) gives
\be\label{papai}
\int \E_\eta (\psi(x,\eta(t))\psi(y, \eta(t)))\
\nu_\rho(d\eta)= (1-\rho)^2\E^{\SEP}_{x,y}\left( e^{-\lambda
\int_0^t  (\delta_{x(s),0}+ \delta_{y(s),0})\ ds }\right).
\ee
Taking the limit $t\to\infty$ gives
\[
\lim_{t\to\infty}\int \E_\eta (\psi(x,\eta(t))\psi(y, \eta(t)))\
\nu_\rho(d\eta)=(1-\rho)^2\E^{\SEP}_{x,y}
\left( e^{-\lambda\int_0^\infty  (\delta_{x(s),0}+ \delta_{y(s),0})\ ds }\right),
\]
which proves 2b) of part I of the theorem.
Item 2b) of part II is proved along the same lines.
\epr
Using the  alternative  duality relation \eqref{fedual},
we can prove the  following  more general statement.
\bc\label{gencorrel}
Let $(p(x,y): x,y\in V)$ be transient, then
for $\mu_{\rho,0}=\lim_{t\to\infty} \nu_\rho S_0(t)$  we have the following.
For all $x_1, \ldots, x_n$, $n$ distinct points in $V$, we have
\beq
&&\int\prod_{i=1}^n (1-\eta_{x_i}) \mu_{\rho,0} (d\eta)- \prod_{i=1}^n\int (1-\eta_{x_i}) \mu_{\rho,0} (d\eta)
\nonumber\\
&=&
(1-\rho)^n\left(\E^{\SEP}_{x_1, \ldots, x_n} e^{-\lambda\sum_{i=1}^n \int_0^\infty\delta_{x_i(s),0}ds}
-\E^{\IRW}_{x_1, \ldots, x_n} e^{-\lambda\sum_{i=1}^n
\int_0^\infty\delta_{X_{i,s},0} ds}\right).\nonumber\\
&&
\eeq
\ec

\section{Further properties of the invariant measures}\label{sec:further}
\subsection{Negative correlations}\label{subsec:neg-cor}
 To derive more properties of the invariant measure $\mu_{\rho,0}$,
we compare the evolutions of exclusion process and of independent random walkers.
We first prove the following lemma,
which will imply negative correlations in the measure $\mu_{\rho,0}$.
Recall the notation introduced in Subsection \ref{subsec:exclusion}.

\bl\label{corlemlig}
For all $t\geq 0$, and for all $x_1, \ldots, x_n$, $n$ distinct points in $V$
\be
\left(\E^{\SEP}_{x_1, \ldots, x_n} e^{-\lambda\sum_{i=1}^n
\int_0^t\delta_{x_i(s),0}ds}
-\E^{\IRW}_{x_1, \ldots, x_n} e^{-\lambda\sum_{i=1}^n
\int_0^t\delta_{X_{i,s},0}ds}\right)
\leq 0.
\ee
\el
\bpr
Let us
denote $\Psi(x_1, \ldots, x_n)= \sum_{i=1}^n \delta_{x_i, 0}$.
Then we have, using the Feynman-Kac formula
\beq\label{bopea}
&&\left(\E^{\SEP}_{x_1, \ldots, x_n}
e^{-\lambda\sum_{i=1}^n
\int_0^t\delta_{x_i(s),0}ds}-\E^{\IRW}_{x_1, \ldots, x_n}
e^{-\lambda\sum_{i=1}^n \int_0^t\delta_{X_{i,s},0}ds}\right)
\nonumber\\
&=&
(e^{t(\caV_n +\Psi)}1- e^{t(\caU_n +\Psi)}1)(x_1, \ldots, x_n)
\nonumber\\
&=&
\int_0^t ( e^{(t-s)(\caV_n +\Psi)}(\caV_n-\caU_n)
e^{s(\caU_n +\Psi)}1)(x_1, \ldots, x_n)ds.
\eeq
 Here in the last equality we used partial integration,
as in \cite{ligg}, Chapter VIII, proof of Proposition 1.7.
Because the function $1$ is positive definite, and the semigroup
$(e^{s(\caU_n +\Psi)}1)(x_1, \ldots, x_n)$ factorizes, i.e.,
\[
(e^{s(\caU_n +\Psi)}1)(x_1, \ldots, x_n)
= \sum_{z_1, \ldots, z_n} \prod_{i=1}^nk_s(x_i, z_i),
\]
where
\[
k_s(u,v)= \E_u^{\IRW} \left(e^{-\lambda\int_0^t \delta_{X_{s},0} ds}I(X_t=v)\right),
\]
we have that $(e^{s(\caU_n +\Psi)}1)(x_1, \ldots, x_n)$
is a positive definite symmetric function of $x_1,\ldots, x_n$.
Therefore, by Liggett's inequality
(see \cite[Chapter VIII, Proposition 1.7]{ligg}),
we have $((\caV_n-\caU_n) e^{s(\caU_n +\Psi)}1)(x_1, \ldots, x_n)\leq 0$.
The result then follows from  \eqref{bopea} using that
$e^{(t-s)(\caV_n +\Psi)}$ is a positive semigroup (i.e., maps
non-negative functions to non-negative functions).
\epr

Then we have
the following corollary  for the process with a source.
\bp[Negative correlations]
Let $(p(x,y): x,y\in V)$ be transient, then
for $\mu_{\rho,0}=\lim_{t\to\infty} \nu_\rho S_0(t)$  we have the following.
For all $x_1, \ldots, x_n$, $n$ distinct points in $V$, we have
\beq
\int\prod_{i=1}^n (1-\eta_{x_i}) \mu_{\rho,0} (d\eta)
- \prod_{i=1}^n\int (1-\eta_{x_i}) \mu_{\rho,0} (d\eta)\leq 0.
\eeq
\ep
\bpr
This follows by combining Corollary \ref{gencorrel} with Lemma \ref{corlemlig}.
\epr

To complement this result, we show that in general  $\mu_{\rho,0}$
is not a product measure.
\bp\label{nonprodprop}
Let $(p(x,y): x,y\in V)$ be transient, then for $0<\rho<1$,
$\mu_{\rho,0}=\lim_{t\to\infty}  \nu_\rho S_0(t)$  is not a product measure.
\ep
\bpr
We have
\[
\E_{\nu_\rho} (1-\eta_x(t))
= (1-\rho) \E^{\RW}_x \left(e^{-\lambda \int_0^t \delta_{X_s,0}\ ds}\right)=: h(x,t).
\]
Then $h(x,t)$ satisfies
\[
\frac{d h(x,t)}{dt}= L_{\RW} h(x,t) -\lambda \delta_{x,0} h(x,t),
\]
with $h(x,0)= 1-\rho$. As a consequence,
$u(x)=\lim_{t\to\infty} h(x,t)$  exists because
$(p(x,y): x,y\in V)$ is transient, and it  satisfies
\be\label{harm}
L_{\RW} u(x)= \lambda \delta_{x,0} u(x).
\ee
Let us denote by $\La$ the product measure with $\int (1-\eta_x) \La(d\eta) = u(x)$.
 If this measure were invariant then, for all functions in the domain of the generator
$L_{\SEP,0}$ we would have $\int L_{\SEP,0} f d\Lambda =0$.

Now fix $x\not=0$ such that $p(0,x)>0$ and compute, using the generator
\eqref{sourcexgen} and \eqref{popi}.
\begin{eqnarray*}
&&\int ( L_{\SEP,0}  (1-\eta_x)(1-\eta_0)) \La(d\eta)
\nonumber\\ &=& (L_{\RW} u(0)) u(x) + (L_{\RW} u(x)) u(0)
 - p(0,x) (u(0)- u(x))^2  - \lambda u(0)u(x)
\nonumber\\
&=& - p(0,x)  (u(0)- u(x))^2 \not=0.
\end{eqnarray*}
Here in the last step we used \eqref{harm}.
So we conclude that $\La$ is not invariant.
Because we proved earlier that $\mu_{\rho,0}$ is invariant,
it cannot be equal to $\La$.
\epr
\br
In \cite{sch1} it is proved that the two-point function on $\Z$ with a source 
and a sink exhibits negative correlations for all finite times. With essentially 
the same proof, Proposition 4.1 holds also when starting from a product measure 
for the measure at any finite time. In this sense, it can be viewed as an extension 
of \cite{sch1}.
\er

\subsection{The covariance}\label{subsec:cov}
To understand better the covariance in the limiting measure $\mu_{\rho,0}$,
we approximate
$\left(\E_{x,y}^{\SEP} e^{-\lambda
 \int_0^\infty  (\delta_{x(s),0}+\delta_{y(s),0})\ ds}
 -\E_{x,y}^{\IRW}
 e^{-\lambda \int_0^\infty (\delta_{X_s,0}+\delta_{Y_s,0})\ ds}\right)$
around $\lambda=0$ up to second order in $\lambda$.
 Let us denote by
 $\eta=e_x+e_y$  the configuration with one particle at $x$ 
 and one at $y$, so that
$ \delta_{x(s),0}+\delta_{y(s),0}$ is the occupation at zero
at time $s$ starting from $e_x+e_y$ initially.
First notice that 
the expression
$\int_0^\infty (\delta_{x(s),0}+ \delta_{y(s),0}) ds$ can be rewritten as
$\int_0^\infty \eta_0(s)ds$ where $\eta_0(s)$ denotes
the number of particles at $0$ at time $s$.\\
The zero-th order of the approximation is clearly zero,
and the first order is zero because we can compute
\beq
\E_{x,y}^{\SEP}( \delta_{x(s),0}+\delta_{y(s),0}) 
&=&
\E_{\xi}^{\SEP}( \eta_0(s))
\nonumber\\
&=&
\sum_{y} p_t(0,y) \eta_y(0)
\nonumber\\
&=&
\E_{\eta}^{\IRW}( \eta_0(s)),\label{firstorder}
\eeq
where in the last two lines we used self-duality for SEP and 
independent random walkers.

In order  to prepare the computation of the second order term, 
let us denote by
$
p^{\SEP}_t( e_x+e_y, e_u+e_v)
$
the transition probability  for two exclusion particles
to arrive at time $t$ at the configuration $e_u+e_v$
when initially started from $e_x+e_y$.
and   $p^{\IRW}_t( e_x+e_y, e_u+e_v)$
for the corresponding independent particles.
Let us restrict to the transient case, where the associated
Green's functions are well defined by

\beq\label{green}
G^{\SEP} (x,y; u,v)
&= & \int_0^\infty p^{\SEP}_t( e_x+e_y, \delta_u+\delta_v)\  dt
\nonumber\\
G^{\IRW}  (x,y; u,v) &=&\int_0^\infty
p^{\IRW}_t( e_x+e_y, e_u+e_v)\ dt
\nonumber\\
G(x,y) &=& \int_0^\infty p_t(x,y)\ dt.
\eeq
We then have the following.
\bp\label{prop:cov}
Assume that the single particle random walk is transient.
In the notation of \eqref{covmu},
for the symmetric exclusion process with source at the origin, 
we have, as $\lambda\to 0$,
\beq
\cov_{\mu_{\rho,0}} (\eta_x, \eta_y)
=\lambda^2(1-\rho)^2 \psi(x,y) + o(\lambda^2),
\eeq
where
\beq
\psi(x,y)&=& 2\sum_{z\not=0 }   G(0,z)
\left( G^{\SEP} (x,y; 0,z) - G^{\IRW} (x,y; 0,z)\right)
\nonumber\\
&-& 2 G(0,0)G^{\IRW}(0,0;x,y).
\eeq
\ep
\bpr
By expanding the exponential
$e^{-\lambda  \int_0^\infty (\delta_{x(s),0}+ \delta_{y(s),0})ds}$
up to second order in $\lambda$, we observe that we have to compute
\beq
&& \int_0^\infty\int_0^\infty
\E^{\SEP}_{e_x+e_y}\eta_0(s) \eta_0(r) \ ds dr
- \int_0^\infty\int_0^\infty \E^{\IRW}_{e_x+e_y}\eta_0(s) \eta_0(r)\ ds dr
\nonumber\\
&=&
2\int_0^\infty\int_0^\infty \E^{\SEP}_{e_x+e_y}\eta_0(s) \eta_0(s+r) \ ds dr 
\nonumber\\
&&- 2\int_0^\infty\int_0^\infty \E^{\IRW}_{e_x+e_y}\eta_0 (s) \eta_0 (s+r)\ ds dr.
\nonumber
\eeq
Here we used the elementary computation
\begin{eqnarray*}
\int_0^\infty ds \int_0^\infty dr f(s) f(r) &= & 2\int_{s>r}
f(s) f(r) ds dr
\\
&=& 2\int_0^\infty dr\int_r^\infty ds f(r) f(s)
\\
&=& 2\int_0^\infty dr\int_0^\infty  dv f(r) f(r+v).
\end{eqnarray*}
For the symmetric exclusion process,
by self-duality (see \eqref{sdualex}), and using also that
$\eta_0(s)^2=\eta_0(s)$ because $\eta_0(s)\in \{0,1\}$,
we obtain the following.
\beq\label{sepdu}
&&\E^{\SEP}_{e_x+e_y}\eta_0 (s) \eta_0({s+r})
\nonumber\\
 &=& \sum_{z\in V, z\not= 0} p_r(0,z) \E_{e_x+e_y}^{\SEP}(\eta_0(s)\eta_z(s))
+  p_r(0,0) \E_{e_x+e_y}^{\SEP}(\eta_0(s))\nonumber\\
&=& \sum_{z\in V, z\not= 0}
p_r(0,z) p^{\SEP}_s(e_x+e_y; e_0+e_z)
\nonumber\\
& +& p_r(0,0)(p_s(0,x)+ p_s(0,y)).
\eeq
By self-duality of independent random walkers (see Subsection \ref{subsec:duality}) we have
\beq\label{irwdu}
&& \E^{\IRW}_{e_x+e_y}(\eta_0(s) \eta_0(s+r))
\nonumber\\
&=&
\sum_{z\in V, z\not=0} p_r(0,z)
\E_{e_x+e_y}^{\IRW} ( \eta_0(s)\eta_z(s))
+ p_r(0,0) \E_{e_x+e_y}^{\IRW} ( \eta_0(s)\eta_0(s))
\nonumber\\
&=&
\sum_{z\in V, z\not=0} p_r(0,z)p_s^{\IRW}(e_x+e_y; 0,z)
+ p_r(0,0)\E_{e_x+e_y}^{\IRW} ( \eta_0(s)(\eta_0(s)-1))
\nonumber\\
&+& p_r(0,0)(p_s(0,x)+p_s(0,y))
\nonumber\\
&=&
\sum_{z\in V, z\not=0} p_r(0,z)p_s^{\IRW}(x,y; 0,z) + p_r(0,0)(p^{\IRW}_s(0,0;x,y))
\nonumber\\
&+ & p_r(0,0)(p_s(0,x)+p_s(0,y)).
\eeq
Subtracting \eqref{irwdu} from \eqref{sepdu} gives
\beq\label{barankol}
&&\int_0^\infty\int_0^\infty
\E^{\SEP}_{e_x+e_y}\left(\eta_0(s) \eta_0(r) \right) \ ds dr
- \int_0^\infty\int_0^\infty \E^{\IRW}_{e_x+e_y}
\left(\eta_0(s) \eta_0(r)\right)\ ds dr\nonumber\\
&=&
\sum_{z\in V, z\not=0} G(0,z)\left(G^{\SEP}(x,y; 0,z)
- G^{\IRW}(x,y; 0,z)\right)
\nonumber\\
&&- \, G(0,0)G^{\IRW}(0,0;x,y).
\eeq
\epr
\br
Proposition \ref{prop:cov}
shows that the correlations in the limiting measure have long-range character, 
because $ G(0,0)G^{\IRW}(0,0;x,y)$ decays as a power law when $|x|,|y|\to\infty$,
which suggests that the stationary measure $\mu_{\rho,0}$ has properties
of a self-organized critical state. By this we mean that the long-range 
correlations are not a consequence of tuning of parameter such as the 
temperature, but they appear spontaneously, because the unperturbed system 
has a conserved quantity. Such long range correlations are also studied 
in perturbations of the exclusion process in \cite{maesredig} (including 
bounded perturbations), and in a translation invariant setting in \cite{maes}.
In non-equilibrium steady states, long-range correlations appear in 
the exclusion process as was first observed in \cite{spohn}.
\er

\section{Independent random walkers}\label{sec:irw}
In this section, we consider the same problem of adding a source
in the context of independent random walkers.
When initially started from a homogeneous Poisson measure,
we show that at any later time $t>0$, we still have a Poisson measure,
and depending on the transience/recurrence of the underlying random walkers,
we obtain a limiting non-homogeneous density profile, or the density
(expected number of particles  at each site) blows up in the limit
$t\to\infty$. This result is proved via an intertwining relation with
a deterministic process, which we believe is of independent interest.
 Notice that the conservation of Poisson product measures is due to the fact 
 that particles are added at constant rate $\lambda$. If we remove particles 
 at rate $\mu$ (sink), then Poisson product measure are no longer conserved, 
 so we will not consider this case here.
\subsection{Independent random walkers: generator}\label{subsec:irw-gen}
As in the previous section, let $V$ denote a countable set of vertices and
$(p(x,y):x,y\in V)$ an irreducible and symmetric random walk transition
probability on $V$.
Recall that  the independent random walkers process is
 a process on $\N^V$ with formal generator \eqref{ind}.
The independent random walkers process with source at $0\in V$ is then
defined via the formal generator
\beq\label{sourceind}
 L_{\IRW,0}  f(\eta) &=&
\sum_{x,y\in V} p(x,y)\left( \eta_x (f(\eta^{x,y})- f(\eta))
 + \eta_y (f(\eta^{y,x})- f(\eta))\right)
\nonumber\\
&+&
\lambda (f(\eta+ e_0)-f(\eta)).
\eeq
The generators \eqref{ind} and \eqref{sourceind} can be rewritten
in terms of creation and annihilation operators defined as follows.
For a function $f:\N\to\R$ we define
\be\label{def:a-adagger}
 af(n)= n f(n-1) I(n\geq 1),\qquad a^\dagger f(n) = f(n+1)\ee
 See also \cite{doi} where these operators were introduced
 in the context of a ``quantum field theory'' for reaction-diffusion systems.
For a local function $f: \N^V\to\R$ we then define the operators
$a_i, a_i^\dagger $ for $i\in V$ as
\beq\label{anni}
a_i f(\eta) &=& \eta_i f(\eta-e_i)\nonumber\\
a_i^\dagger f(\eta) &=& f(\eta+e_i).
\eeq
where $\eta_i f(\eta-e_i)$ is a shorthand for
$a_i f(\eta)= 0$ if $\eta_i=0$ and
$a_i f(\eta)= \eta_i f(\eta-e_i)$ otherwise.
The generator \eqref{sourceind} can then be rewritten as follows
\be\label{sourceann}
L_{\IRW,0}= -\sum_{x,y\in V} p(x,y) (a_x- a_y)(a_x^\dagger- a_y^\dagger)
 + \lambda (a_0^\dagger- Id),
\ee
where $Id$ denotes the identity operator.
\subsection{Intertwining operators}\label{subse:interop}
In this subsection we prove that the process with generator
\eqref{sourceind} starting from $\nu_\rho$  is intertwined with 
a deterministic process.
This amounts to find a differential operator representation
for the creation and annihilation operators \eqref{anni}, in the spirit of
\cite{gkrv}.
In order to define the intertwiner, we first introduce for a function
$f: \N\to\R$ the associated generating function
\be\label{bolanko}
Gf(z)= \sum_{n=0}^\infty f(n) \frac{z^n}{n!} e^{-z},
\ee
where we implicitly assume that the function $f$ is such that
the series defining $Gf(z)$ is absolutely convergent in an open interval
around the origin.
Notice that for $z\geq 0$, $Gf(z)$ is precisely the expectation of $f$ w.r.t.
Poisson distribution with parameter $z$.
The following proposition collects intertwining relations between $G$
and the creation and annihilation operators defined above.
\bp\label{intprop}
Let $a, a^\dagger$ be defined as above in \eqref{def:a-adagger},
and $G$ the generating function \eqref{bolanko}. Then we have
\beq
G (a f) (z) &= & z Gf(z)
\nonumber\\
G (a^\dagger f) (z) &=& \frac{\partial Gf}{\partial z} (z) + Gf(z).
\eeq
\ep
\bpr
This follows from a simple computation.
\epr

\noindent
We now extend the generating function $G$ to the multi-variate setting
by  tensorization.
  Let $\widehat\Omega$ denote
the configuration space $\R^V$. With small abuse of notation
we denote by $z$ configurations in $\widehat\Omega$. 
For $f:\N^V\to\R$ a local function depending on $\eta_i, i\in \la$
with $\la\subset V$ finite, we define, 
\be\label{gee}
\caG f(z)  = \sum_{\eta_i, i\in \la} \prod_{i\in\la}\frac{z_i^{\eta_i}}{\eta_i !}
e^{-z_i} f(\eta_i, i\in \la),
\ee
 where $z=(z_i,i\in V)\in \widehat\Omega$. 
Notice that if $z_i\geq 0$ for all $i\in V$, 
then $\caG f(z)= \nu_z (f)$ where, again with a small abuse of notation,
$\nu_z=\otimes_{i\in V} \nu_{z_i}$ is the product
of Poisson measures with parameter $z_i$ at site $i\in V$.

We have the following intertwining relation.
Denote, for $f$ local and smooth, and $z\in \widehat\Omega$
\be\label{interz}
\loc f(z)= \sum_{x,y\in V} - p(x, y) (z_x- z_y)(\partial_x-\partial_y) f(z)
 + \lambda (\partial_0 f)(z),
\ee
with $\partial_x= \frac{\partial}{\partial z_x}$.
Then $\loc$ is the generator of a deterministic system of differential equations
\be\label{bakok}
\frac{dz_x(t)}{dt}= \sum_{y\in V} p(x,y) (z_y(t)- z_x(t)) + \lambda \delta_{x,0}.
\ee
Using the generator \eqref{Agen} and the associated semi-group,
we can rewrite \eqref{bakok} as follows
\be
\frac{dz_x(t)}{dt} =  (L_{\RW} z)_x(t) + \lambda \delta_{x,0}.
\ee
This equation can be solved by the classical variation
of constants method, and we obtain
\be\label{bakaro}
z_x(t)= \E_x (z_{X_t}(0)) + \lambda\int_0^t \E_x (\delta_{X_s, 0} ) ds.
\ee
Using proposition \ref{intprop} we then obtain the following
intertwining relation and as a consequence evolution of Poisson product measures.
\bp\label{intpropo}
Let  $L_{\IRW,0}$  denote the generator of the independent random walkers with source
at the origin given in \eqref{sourceind}  and
$\loc$ the generator of the deterministic system given in \eqref{interz}.
Then we have, for every local function $f: \N^V\to\R$, such that
$\caG f$ exists and is smoothly depending on $z$:
\be\label{Lloc}
\caG L_{\IRW,0} f = \loc \caG f.
\ee
As a consequence for any $\rho: V\to [0,\infty)$, if we denote $\nu_\rho$
the product of Poisson measures
with parameter $\rho(x)$ at $ x\in V$, and $\nu_\rho (t)$ the measure obtained
at time $t$ when starting
the Markov process with generator \eqref{sourceind} 
 and initial distribution $\nu_\rho$,  then we have
\be\label{lasteq}
\nu_{\rho} (t) = \nu_{\rho_t},
\ee
with
\be\label{lastsol}
\rho_t(x)= \E_x (\rho(X_t)) + \int_0^t \E_x (\delta_{X_s, 0} ) ds,
\ee
the solution of \eqref{bakaro} with initial condition
$z_x(0)= \rho(x)$.
The density $\rho_t(x)$ diverges as $t\to\infty$ if and only if
the random walk with generator  $L_{\RW}$  is recurrent.
If this random walk is transient, starting from a homogeneous Poisson
product measure with density $\rho$, $\nu_\rho(t)$ converges, as
$ t\to\infty$ to the Poisson product measure with density
\[
\rho_\infty(x)= \rho + \lambda\int_0^\infty \E_x (\delta_{X_s, 0} ) ds,
\]
which is a solution of the equation
\[
L_{\RW} f(x)= \lambda \delta_{x,0}.
\]
\ep
\br
Connections between deterministic systems of linear differential equations and
a system of independent random walks were studied in \cite{pelity}, see also
\cite{gkrv} where a duality between a system of independent random walkers and a deterministic system is proved. Here, on the contrary, we prove an intertwining relation, which can be viewed as a way of recasting and extending Doob's theorem on the conservation of Poisson product measures under independent evolutions, see e.g.\ \cite{dmp}.
\er
\noindent
{\bf Acknowledgements.}
We thank the referees for their careful reading and useful feedback.
We thank Mario Ayala, Jochem Hogendijk and Mitchell Maassen Van den Brink
for useful discussions.  F.R. thanks MAP5 lab. for financial support
and hospitality.

\end{document}